# Tromino Tiling Deficient Cubes of Any Side Length

Norton Starr (Amherst College, Amherst, MA, USA)

ABSTRACT. We show that three dimensional cubes of any size can be tiled with trominoes and, when necessary, one or two singletons in any positions. Cubes of side length a multiple of three can always be tiled with trominoes (known), cubes of side length congruent to 1 mod 3 can always be tiled with an arbitrary single cube and trominoes, and cubes of side length congruent to 2 mod 3 can always be tiled with two single cubes in arbitrary locations and trominoes.

Solomon Golomb gave an elegant proof by induction for tiling a class of checkerboards by one type of polyomino: He showed that any deficient n×n square grid having side length a power of 2 can be tiled by plane trominoes [**4**]. A (plane) tromino consists of three congruent squares of unit area, two of which meet the third on adjacent edges as shown on the left in Figure 1:

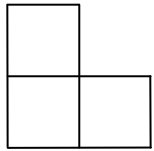 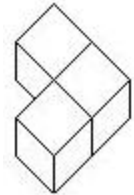

**Figure 1**

An n×n square is *deficient* if one of its cells is occupied or otherwise removed from consideration. To *tile* a deficient square having side length n is to fill the remaining $n^2-1$ cells by trominoes.

Chu and Johnsonbaugh [**3**], Ash and Golomb [**1**] and others have studied extensions of this problem to arbitrary squares and to rectangles with two deficiencies as well. The key to Golomb's mathematical induction is concisely demonstrated by Roger Nelsen's "proof without words" in Figure 2*. Using perpendicular slices through the center to partition the deficient square into four quarters, the second grid shows how the appropriate insertion of a tromino about the center reduces the tiling of the original square to the tiling of four deficient squares, each having side length half that of the first. The case of an 8×8 square is now a popular recreation sold by Kadon Enterprises under the name Vee-21. The history of this 8×8 puzzle is to appear in [**7**].

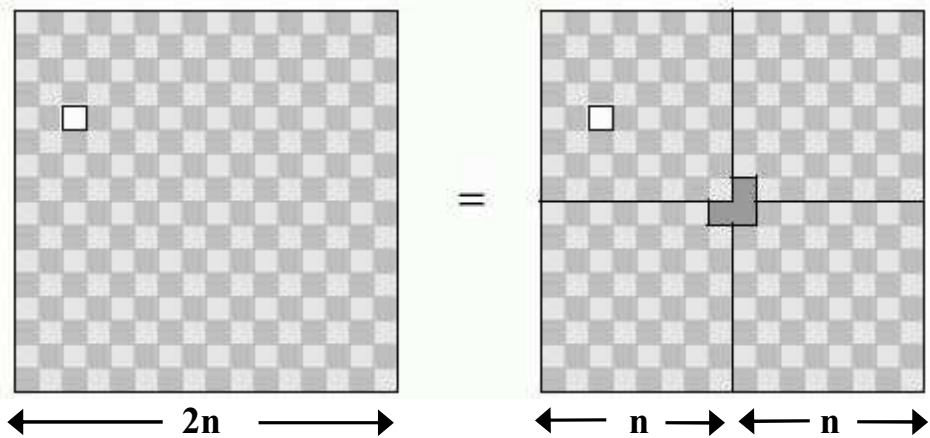

**Figure 2***

*From *Proof Without Words II* [**6**], with the permission of the Mathematical Association of America and Roger B. Nelsen. © 2000 The Mathematical Association of America.

Our aim is to generalize this result to three dimensions. Take as basic building blocks cubes of side length one (called *singletons*) and clusters of three such cubes (called *trominoes*) consisting of a singleton together with two other singletons meeting it on adjacent faces, as shown on the right in Figure 1. Three dimensional n×n×n cubes are *singly deficient* if one of their $n^3$ cells is occupied or otherwise removed from consideration, and *doubly deficient* if two of the cells are occupied or otherwise removed from consideration. The goal is to fill the remaining cells by trominoes, namely, to *tile* the cube.

For distribution at a puzzle conference in 2007 (IPP27) we designed a singly deficient 4×4×4 puzzle, to be tiled with one singleton and 21 trominoes. This is a simple exercise and proving that it can be carried out for any placement of the singleton is also easy. The corresponding problem of tiling deficient 5×5×5 cubes forces at least two singletons, because $5^3 \equiv 2 \pmod 3$. Proving that such a cube can be tiled was a daunting task, because the number of arbitrary placements of two deficiencies seemed too many to merit study. However, doubling the side length from 4 to 8 suggested the possibility of a recursive tiling argument akin to that of Golomb's method illustrated in Figure 2 for deficient squares. Though $8^3$, like $5^3$, is congruent to 2 mod 3, a proof using recursion to cubes of side length 4 is indeed feasible, if somewhat tedious. With the results for side lengths 4 and 8 in hand, we proved that any singly deficient cube of side length $2^{2n}$ can be tiled with trominoes, and any doubly deficient cube of side length $2^{2n+1}$ can be tiled with trominoes. The proof used the Golomb method together with a variety of strong induction. For even exponents, the recursive argument was easy, while for odd exponents it was complicated by the extra degree of freedom a second singleton provides. (Note that as a consequence, any singly deficient cube of side length $4^n$ can be tiled by trominoes.) These arguments are given in [**8**].

The result for cubes of side length $2^k$ raised the obvious question of tiling cubes of arbitrary side length n. It was easy to show that cubes having side length a multiple of 3 could be tiled using only trominoes – no singletons were necessary. Indeed, this had already been presented by Boltyanski and Soifer as a special case of Exercise 6.10 on page 59 of [**2**]. Experimenting with 41 trominoes and two singletons suggested that all doubly deficient cubes of side length 5 can be tiled. Observing that 5 (as well as 8) is congruent to 2 mod 3, 4 is congruent to 1 mod 3, and 3 is congruent to 0 mod 3 suggested that the residue mod 3 determines the number of singletons necessary and sufficient for a tromino tiling. Our main theorem will express this: Cubes of side length a multiple of three can be tiled with trominoes, all singly deficient cubes of side length congruent to 1 mod 3 can be tiled with trominoes, and all doubly deficient cubes of side length congruent to 2 mod 3 can be tiled with trominoes.

We give the proof in several stages. Cubes of side length n =1 and 2 are trivial (the latter is shown in Lemma 1 of [**8**]) and the easy case n = 4 is Theorem 1 of [**8**]. For side lengths congruent to 0 mod 3 the argument is easy, and appears in Theorem 2 below. The messy case n = 5 is given below as Theorem 1. Next we show how to tile some basic rectangular solids that will be used repeatedly in the induction arguments proving the General Theorem. With base cases n = 3, 4, and 5, we then use strong induction, explaining how to recursively relate the tiling of a cube of side length n + 3 to that of a cube of side length n. This we do in stages: first no singleton (Theorem 2), then singly deficient cubes (Theorem 3), and finally doubly deficient cubes (Theorem 4). The most complicated arguments are for n = 5 and for doubly deficient cubes in general.

**Theorem 1.** Any doubly deficient cube 𝔇 of side length 5 can be tiled with trominoes.

*Proof.* Consider a 4×4×4 sub-cube ℭ contained in 𝔇 and having a corner in common with 𝔇. Of the eight such positions for ℭ, at least one contains at least one of the singletons **S**, **T**.

## Case A

Suppose singleton **S** lies in the 4×4×4 sub-cube ℭ and **T** is outside ℭ. By Theorem 1 of [8], ℭ can be tiled by trominoes. We may assume ℭ lies in the lower, left rear of 𝔇 as we view 𝔇 from the front. The 61 cells in 𝔇 that lie outside ℭ form a shell one unit thick, consisting of three adjacent faces of the cube 𝔇. The first image in Figure 3 is an overhead view of this shell, which may be regarded as a disjoint union of a 4×4×1 top, a 5×4×1 right face, and a 5×5×1 front face. The large 4×4 square in this first image is the part of the shell lying atop ℭ. The 4×1 rectangle on the right in this first image is an overhead view of the 5×4×1 right face. The 1×5 rectangle along the bottom is an overhead view of the 5×5×1 front face. **T** must lie in one of these three faces and we may assume it's in the front 5×5 face: If **T** lies in the right face or in the top face, a tiling follows by appropriate rotations of a tiling for the case where it's in the front face. We now show how to tile the shell with 20 trominoes, no matter where **T** is located in the front face.

On the right of Figure 3 is a numbered grid for reference to particular cells in the front face.

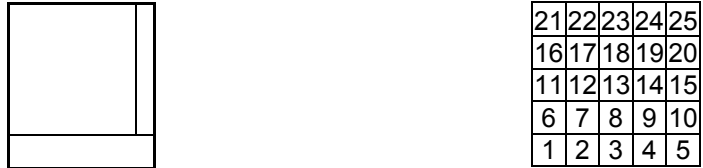

**Figure 3**

Figure 4 shows three tilings by trominoes of the front 5×5 face, for the cases where **T** is in the center of an edge, in the center of the face, or in a corner of the face. (Cf. [**5**], pp. 28, 29.) These are the only positions for **T** admitting a tiling of the face, as shown by Martin ([**5**], p. 28) and by Ash and Golomb ([**1**], p. 51).

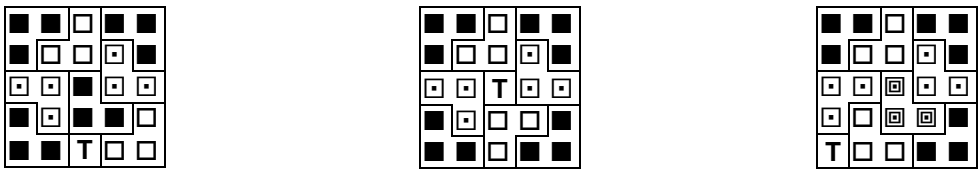

**Figure 4**

For any of such tiling of the front face, it's easy to tile the remaining two faces. Figure 5 shows the top face as a 4×4 grid, viewed overhead *from the right* of 𝔇, and the right face as a 5×4 grid seen also from the right of 𝔇. Neither face can be tiled by trominoes, but with the tromino indicated by **R** straddling these two faces, a complete tiling as shown is possible. Thus, when **T** lies in any of the five positions on the front face as described in the above paragraph, a tiling of 𝔇 is possible: Locations 1, 3, 5, 11, 13, 15, 21, 23 and 25 for T admit tiling.

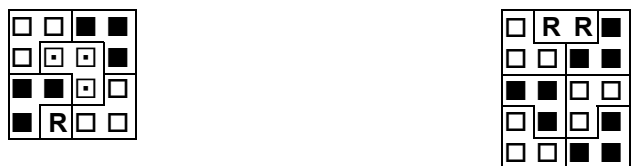

**Figure 5**

Now suppose T lies in one of the 16 cells in the front 5×5 face that prevent tiling of this face. These 16 lie in either the second or fourth row or the second or fourth column or both. We tile this face by having some trominoes protrude into the top or the right face. It's possible to reduce the tiling of the 16 remaining cases to two basic types of protrusions of trominoes into another face, from which any other necessary protrusions can be obtained by reflection. After showing how to tile the front 5×5 face in these two basic cases, we show how to tile the other two faces, completing a tiling of the 5×5 cube and hence all of 𝔇.

The first grid in Figure 6 shows **T** in cell 12, a partial tiling of the front face by trominoes, and three cubes (indicated by **P**'s and **Q**) from two trominoes whose other parts protrude into the 4×4 top face. The second grid shows a similar tiling when **T** is in cell 8, though now the two trominoes partially visible (the **P**'s and **Q**) protrude into the 5×4 right face. We indicate beneath each grid the cell numbers for positions of the singleton **T** accounted for by that grid.

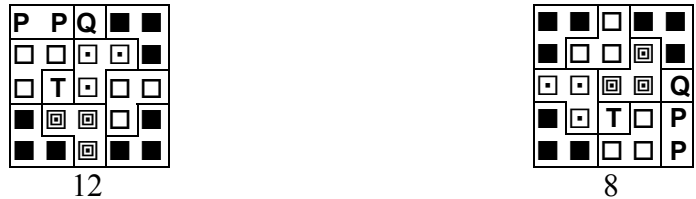

**Figure 6**

For the first grid in Figure 6, we show in Figure 7 how to tile with trominoes the top and the right faces of the shell. This tiling can be reflected to demonstrate a tiling for the second grid of Figure 6: The tiling used on top 4×4 face gets reflected onto the lower 4×4 grid of the right face, and the tiling used on the right 5×4 face gets reflected onto the top of the shell. In both cases the reflections are about the upper right front-to-back edge of the cube 𝔇 as viewed from the front.

The first grid in Figure 7 is an overhead view from the right of the top 4×4 face together with the top edge of the front 5×5 face (the leftmost column in the first grid) and the top edge of the right 5×4 face (beneath the double lines at the bottom of the first grid.) The two trominoes indicated by **P** and **Q** form a 2×3 box, with three of their six cubes lying in the top 4×4 face. Of the 13 remaining cells in the top 4×4 face, 12 are filled by four trominoes and one cube of a fifth tromino (marked by **R**), whose other two cubes lie in the right 5×4 face. The second grid shows the right side's 5×4 face, tiled by 6 trominoes and two of the cubes from the "**R**" tromino. Thus, the entire shell and therefore all of 𝔇 has been tiled: Locations 8 and 12 for **T** admit tiling.

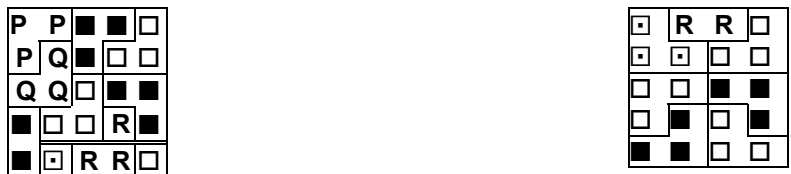

**Figure 7**

Figure 8 shows three views of the front face, each with a column of three cells (15, 20 & 25) occupied by parts (**P** and **Q**) of trominoes that protrude into the right face. The first shows a tiling with **T** in cell 2. Note that if we rotate the 2×2 square containing cells 1, 2, 6 & 7, occupied by **T** and the tromino marked by ▣'s, **T** can be moved to cells 6 or 7. (Ignore cell 1: that case was treated above.) The second grid differs from the first only in having the lower two rows rotated about a vertical line through the middle of the face. Note that the second grid similarly allows for **T** in cells 4, 9 & 10. The third grid shows tiling that allows for **T** in cells 16, 17 & 22.

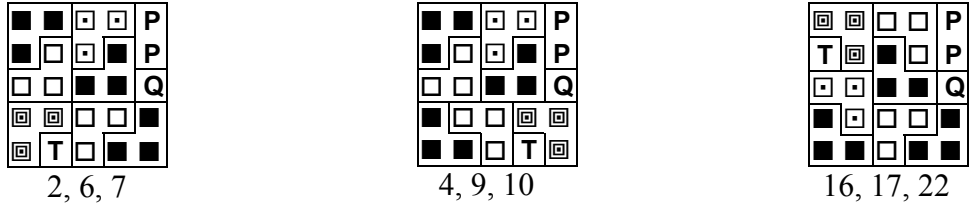

     2, 6, 7                            4, 9, 10                      16, 17, 22
**Figure 8**

Completing the tiling in the above three cases is carried out as shown in the Figure 9, which displays, as viewed from the right, the 4×4 top face and the 5×4 right face. The latter includes three cubes that protrude into it from the front 5×5×5 face, as well as two cubes (marked **R**) whose third cube lies in the 4×4 top face.

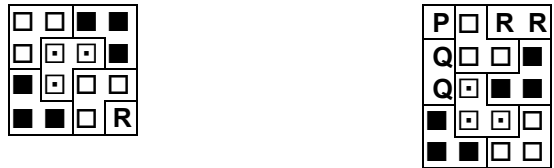

**Figure 9**

Thus locations 2, 4, 6, 7, 9, 10, 16, 17, and 22 for T admit tiling.

     It remains only to show how to tile the front face for T located in cells 14, 18, 19, 20 and 24. We do this by the tilings of the front face given in Figure 10. In the left grid, rotating the 2×2 square occupying cells 18, 19, 23 & 24 (containing **T** and the tromino marked by ▣) lets **T** occupy cells 18, 19 and 24. In the right grid, rotating the 2×2 square occupied by **T** and the ▣ tromino lets **T** occupy cells 14 and 20. The tiling of the remaining two faces (the top and right sides) of the shell for the first grid of Figure 10 is shown in Figure 9 above. For the second grid, the reflection described just below Figure 6 interchanges the tilings of the top and right faces for the first grid to provide tilings that work for the second grid's configuration.

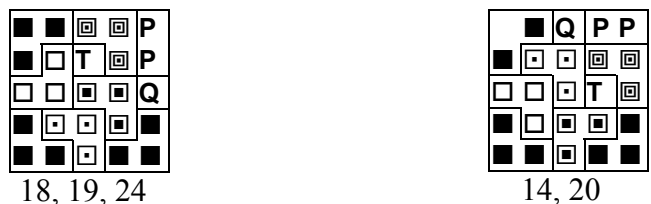

     18, 19, 24                      14, 20
**Figure 10**

     Together with these last five locations: 14, 18, 19, 20 and 24 for **T**, we've now considered all twenty-five possibilities. Therefore we can tile 𝔻 for any **S** is in ℂ and any **T** in the front 5×5 face of the shell consisting of cubes in 𝔻 but outside ℂ. As previously remarked, for **T** on the top or right face of the shell, a rotation so the face containing **T** becomes the front face shows how to complete the tilings in those two cases.

## Case B

Suppose singletons **S** and **T** both lie within the 4×4×4 sub-cube 𝕮 of 𝕯. Tiling the remaining 62 cells of 𝕮 with trominoes forces one or more trominoes to protrude outside 𝕮. Assume **S** and **T** lie in different *horizontal* levels of 𝕮. (Denote the levels A, B, C, and D, bottom to top.) We show how to position cubes protruding from 𝕮 so they lie in the front 5×5 face of the shell in such a way that the rest of the front 5×5 face of 𝕯 can be tiled with trominoes. Then, as shown in Figure 5, the remaining two faces of the shell can be tiled by trominoes, completing the tiling of 𝕯. Finally we show how to handle the case where **S** and **T** lie in the same horizontal layer of 𝕮.

First note that each the two 4×4 levels of 𝕮 containing a singleton can be tiled by Golomb's method for deficient plane squares. This device was displayed in Figure 2 and is also shown in Figure 11, which demonstrates the tiling of a 4×4 grid containing a typical singleton **S**. The key is to place a tromino (marked by ■'s) around the center so it has one cube in each of the three quadrants not containing **S**. Then it's trivial to complete the tiling (as shown by tromino tiles marked with □, ⊡, ⊞, and ⊟.)

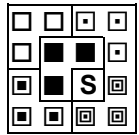

**Figure 11**

Suppose that the two levels containing *neither* **S** *nor* **T** are adjacent. This can happen three ways: These levels could be A and B, or B and C, or C and D. The combined 32 positions in these two levels of 𝕮 can be tiled with a total of 10 trominoes plus two cubes of an eleventh tromino, one of whose cubes protrudes into the front face of the shell. This can be accomplished by first tiling each of these two levels with five trominoes in the exactly same way, so as to leave an empty cell on the front face of 𝕮: An overhead view of this tiling of a level of 𝕮 is given in the left grid of Figure 12. Place a tromino (marked **P**) so two of its unit cubes fill this pair of empty cells on the two adjacent levels, with the third cube protruding into the front 5×5 face of 𝕯. Recall that if this third cube is a singleton lying in the first column of a 5×5 grid, tiling the grid is possible when it's in the 1$^{st}$, 3$^{rd}$ or 5$^{th}$ rows (cells 1, 11 or 21 in Figure 3). Thus we position the **P** tromino so its protruding cube occupies cell 1 (if A, B were the levels containing neither **S** nor **T**) or cell 11 (if B,C or if C, D were the empty levels.) The second grid in Figure 12 shows the front face 5×5 of 𝕯, for the case where levels C, D lacked **S**, **T**. The same grid applies if B, C had been empty, but were A, B empty, the protruding cube would be made to occupy cell 1, which, as shown in the third grid, is also is a 5×5 configuration that's tilable.

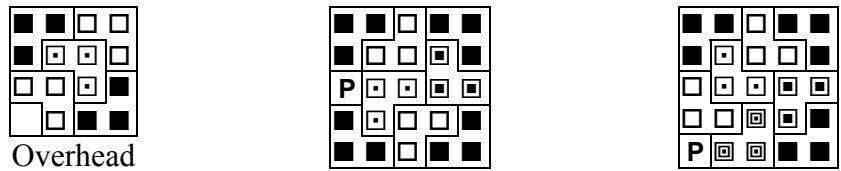

Overhead

**Figure 12**

Suppose that the two levels containing neither **S** nor **T** are A and C: separated by one level. Tile A and C as shown on the left in the above Figure 12, with a single empty cell at the front of $\mathbb{C}$. Using two trominoes, place a single cube from each in these two empty cells, with the other two cubes (marked **P**) in the front 5×5 face of the shell as shown in the first grid of Figure 13. This first grid also demonstrates tromino tiling of the remaining 21 cells on the front 5×5 face. If, instead, B and D are the empty levels separated by one level, tile these in exactly the same way as the previous tiling of A and C, inserting the two **P** trominoes so they have the same protrusions into the front 5×5 face. For the one remaining case, where A and D are the two levels of $\mathbb{C}$ without any singleton, the corresponding protrusions into and tromino tiling of the front 5×5 face are shown in the second grid.

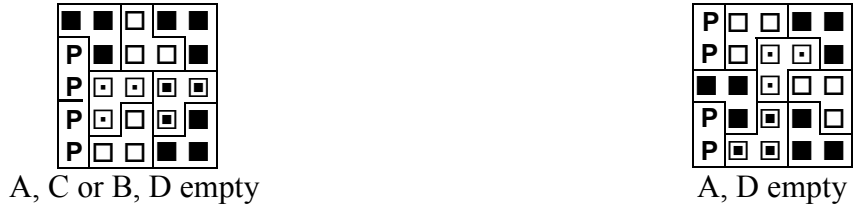

A, C or B, D empty    A, D empty

**Figure 13**

As shown in Figure 5, the above tilings of $\mathbb{C}$ and the front 5×5 face for the cases where **S** and **T** lie on different horizontal levels of $\mathbb{C}$, can be completed to full tilings of $\mathbb{D}$.

The last step in bootstrapping the tiling of $\mathbb{C}$ to that of $\mathbb{D}$ is the case for which **S** and **T** lie in the same horizontal layer of $\mathbb{C}$. A rotation of $\mathbb{D}$ will give an orientation of $\mathbb{C}$ for which, relative to the horizontal, **S** and **T** then lie on different levels. Carrying out the tilings of $\mathbb{C}$ and the new front 5×5 face as done above and then using the patterns of Figure 5, gives the tiling of $\mathbb{D}$. This completes the proof of Theorem 1.

## Cubes of Arbitrary Side Length

**General Theorem** Let n be a positive integer and consider a cube of side length n.
   If n is a multiple of three, the cube can be tiled with trominoes.
   If n is congruent to 1 mod 3, the cube can be tiled with trominoes if it's singly deficient.
   If n is congruent to 2 mod 3, the cube can be tiled with trominoes if it's doubly deficient.

**Remarks** We've noted that cases n = 1, 2 and 4 are dealt with in [**8**], and the case n = 5 was proven in Theorem 1 above. For side lengths congruent to 0 mod 3 the argument is easy, and appears in Theorem 2 below. Our main effort in proving the General Theorem will be the proofs by induction for side lengths congruent to 1 mod 3 (Theorem 3) and congruent to 2 mod 3 (Theorem 4.) The basic structure underlying these recursive proofs is illustrated in Figure 14, which shows the partition of an (n+3)×(n+3)×(n+3) cube into an n×n×n cube, three 3×n×n square "slabs", three 3×3×n "rods", and one 3×3×3 cube.

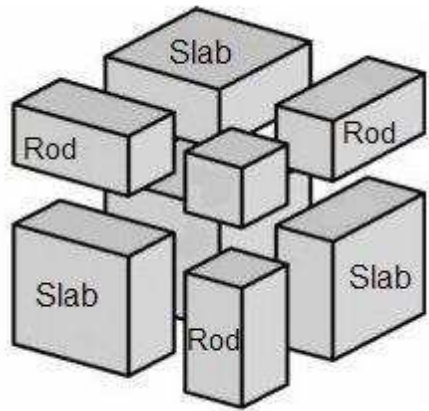

**Figure 14***

*This Figure is a relabeled version of one developed by Thomas Banchoff of Brown University and Davide Cervone of Union College for the April 2000 Math Awareness Month, and appears here with their permission.

We will use tilings of several basic rectangular solids, as given in the following lemmas. They encompass the easy case of cubes having side lengths congruent to 0 mod 3 (Theorem 2 below.). As noted above, this result is a special case of Exercise 6.10 on page 59 of [**2**].

**Lemma 1.** A 3×3×3 cube can be tiled with trominoes.

*Proof.* We demonstrate by pictures a way to accomplish the tiling. The left image in Figure 15 shows how to tile a 3×2×1 box with two trominoes. The middle image in Figure 15 shows how to place three 3×2×1 boxes so as to fill all but nine cells of a 3×3×3 cube: One box is flat and two are lying on their longer edges. The right image shows three parallel trominoes arranged so they can be picked up en masse and inserted among the three 3×2×1 boxes of the middle image to complete the 3×3×3 cube.

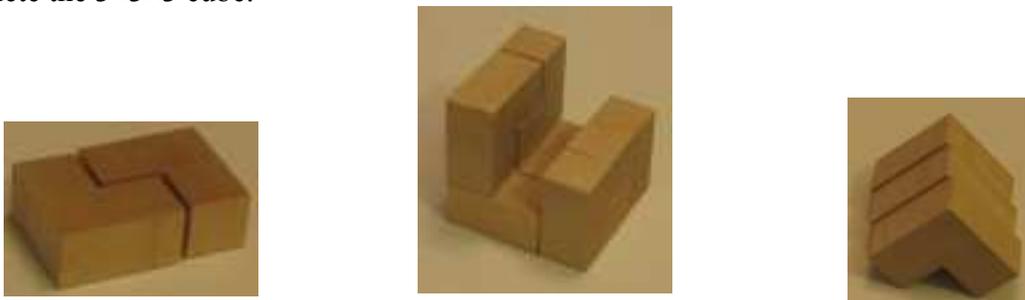

**Figure 15**

**Theorem 2.** For n a positive integer, any cube of side length 3n can be tiled with trominoes.

*Proof.* Such a cube is a union of $n^3$ 3×3×3 cubes.

**Lemma 2.** For n a positive integer, a 3×2×n box can be tiled with trominoes

*Proof.* Place n copies of a 3×2×1 box in parallel, with their largest faces adjacent.

**Lemma 3.** For n a positive integer ≥ 2, the 3×3×n rod can be tiled with trominoes.

*Proof.* By Lemmas 2 and 1, the 3×3×2 and 3×3×3 rods can be tiled. For n even, use n/2 boxes of form 3×3×2. For odd n, using a single 3×3×3 cube suffices if n = 3 and viewing the 27 cubes at one end of the rod as a 3×3×3 cube reduces the rest of the rod's tiling to the case of an even length if n > 3.

**Lemma 4.** For n a positive integer, a 3×(3n+1)×(3n+1) slab and a 3×(3n+2)×(3n+2) slab can be tiled with trominoes.

*Proof.* If n is odd, note that 3n+1 is even, so the 3×(3n+1)×(3n+1) slab is a union of boxes of form 3×2×(3n+1), to each of which Lemma 2 can then be applied.

Suppose n is even. View the 3×(3n+1)×(3n+1) slab as a union of a 3×(3n−2)×(3n−2) slab, two 3×3×(3n−2) rods, and a 3×3×3 cube (as in the top layer of see Figure 14.) Lemmas 1 and 3 assure the tiling of the 3×3×3 cube and the two rods, resp. Because n is even, 3n−2 is even and ≥ 2, so the 3×(3n−2)×(3n−2) slab is a union of boxes of form 3×2×(3n−2), and these can be tiled by Lemma 2.

Essentially the same argument applies to the 3×(3n+2)×(3n+2) slab.

**Theorem 3.** For n a positive integer, any singly deficient cube of side length 3n+1 can be tiled with trominoes.

*Proof.* The result is true for cubes of side length 4 by Theorem 1 of [**8**]. Assume that for a general positive integer n, singly deficient cubes of side length 3n+1 can be tiled with trominoes. Consider a singly deficient cube 𝔇 of side length 3n+4. There are eight ways a cube ℭ of side length 3n+1 can lie within 𝔇 if we require that these two cubes have a corner in common. Choose such a position for ℭ so that it contains the singleton, which we denote T. (At least one position must contain T, because the side length of ℭ is larger than half the side length of 𝔇.)

It may help to envision the structures if we assume 𝔇 has been rotated so that the ℭ containing T is in the lower, left corner at the rear, as suggested in Figure 14. The induction hypothesis says that this lower, left rear (3n+1)×(3n+1)×(3n+1) cube ℭ can be tiled by trominoes. We need to tile the remaining cells of 𝔇. To do this, observe first that the cells of 𝔇 lying outside ℭ form a three-sided "shell" consisting of the following components:

i) A 3×3×3 cube in the upper, right, front corner of 𝔇 (and which meets ℭ at one corner.)

ii) Three square slabs of thickness 3 units and side length 3n+1 (these meet three adjacent, square, (3n+1)×(3n+1) faces of ℭ.)

iii) Three rods of size 3×3×(3n+1) each, which fill the gaps between the three pairs of adjacent square slabs described in ii).

Now tile the shell as follows. Use Lemma 1 to tile the 3×3×3 corner cube, Lemma 3 to tile the rods, and Lemma 4 to tile the (3n+1)×(3n+1)×3 slabs.

Note that the volume of ℂ plus that of the three slabs, together with the three rods and the corner cube, equals the volume of 𝔇: $(3n+1)^3 + 3\cdot(3\cdot(3n+1)^2) + 3\cdot(3^2\cdot(3n+1)) + 3^3 = (27n^3 + 27n^2 + 9n + 1) + (81n^2 + 54n + 9) + 81n + 27 + 27 = 27n^3 + 108n^2 + 144n + 64 = (3n+4)^3$.

We now analyze the most difficult case, the tiling of doubly deficient cubes of side length 3n+2. The proof will follow the lines of the above theorem for singly deficient cubes of side length 3n+1. The hard part occurs when one singleton **S** lies in ℂ and the other, **T**, is in the shell. We will show that if **T** belongs to a 3×2×2 box within a slab or belongs to a 3×3×2 box within a rod or the corner cube, then 𝔇 can be tiled. To do this, we show that the shell can be constructed so that **T** lies in such a box. This is facilitated by the following strengthening of Lemma 3.

**Lemma 5.** (Rod tiling) For n a positive integer, any 3×3×(3n+2) rod can composed of 3×3×2 boxes and at most one 3×3×3 cube.

*Proof.* From the proof of Lemma 3, the rod can be completely tiled using only 3×3×2 boxes if n is even, and using 3×3×2 boxes together with one 3×3×3 cube if n is odd.

**Theorem 4.** For n a positive integer, any doubly deficient cube of side length 3n+2 can be tiled with trominoes.

*Proof.* We use mathematical induction. The base case, n = 1, was demonstrated in Theorem 2, which dealt with cubes of side length 5. Assume the result is true for cubes ℂ of side length 3n+2, and consider a doubly deficient cube 𝔇 of side length 3n+5, containing singleton cubes **S** and **T**. There are eight ways a cube ℂ of side length 3n+2 can lie within 𝔇 if we require that these two cubes have a corner in common. Represent 𝔇 by as shown in Figure 14: It's a union of ℂ, three slabs of form 3×(3n+2)×(3n+2), three rods of form 3×3×(3n+2), and one 3×3×3 corner cube. We now show 𝔇 can be tiled by trominoes.

If there is such a position for ℂ so that it contains both **S** and **T**, then by hypothesis, ℂ can be tiled with trominoes. The shell can also be tiled with trominoes: the slabs can be tiled by Lemma 4, the rods by Lemma 3, and the corner cube by Lemma 1. Thus in this case, 𝔇 can be tiled.

Suppose now that for each position of ℂ, it contains one of **S**, **T**, but not both. We may assume **S** lies in ℂ, so **T** is in 𝔇 but not in ℂ: it's in the shell. We show how to place a tromino **P** having one cube in ℂ and its other two in the shell, in such a way that the rest of the shell can be tiled with **T** and trominoes. Because ℂ then contains **S** and the one cube from the tromino **P**, it is doubly deficient and thus, by the induction hypothesis, it can be tiled, completing the tiling of 𝔇.

There are three possibilities, which we condense to two cases: **T** lies in a 3×(3n+2)×(3n+2) slab, **T** lies in a 3×3×(3n+2) rod, or **T** lies in the 3×3×3 cube. We solve the slab case first, then solve both the rod and cube possibilities by viewing **T** as belonging to an extended rod of shape 3×3×(3n+5). Note that the other six portions of the shell can be tromino tiled by Lemmas 1, 3, and 4.

## Case A. T lies in a (3n+2)×(3n+2)×3 slab

**Let T lie in a (3n+2)×(3n+2)×3 slab and 3n+2 be even**.

Orient 𝔇 so the slab is horizontal, on top of ℭ. Ignoring **T**, this slab can be tiled with 2×2×3 boxes, placed vertically so their 2×2 bases rest atop ℭ, as shown in Figure 16 for the case n = 2:

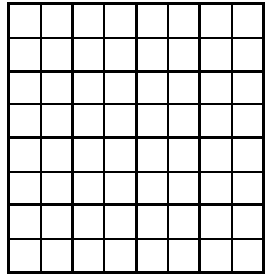

**Figure 16**

Consider the 2×2×3 box B in one of whose 12 cells **T** is located. We now show how to tile the remaining 11 cells of B with three trominoes and two cubes of a fourth tromino whose cubes are marked "**P**", because one its cubes protrudes down into the top layer of ℭ, beneath the bottom 2×2 face of B. Note that the singleton **S** in ℭ cannot already occupy the position in ℭ's face just below B, because if it did, ℭ could have been oriented so as to enclose both **S** and **T**: **S** and **T** would be among the four outer layers of 𝔇, yet ℭ has side length ≥ 5.

Figure 17 shows, from left to right, overhead views of the top, middle and bottom levels of the 2×2×3 box B for the case where **T** lies in the top level of B. The positions of the two cubes of **P** are also displayed, and its third cube would be in ℭ, one level beneath B. In each case, a rotation of B (and the single cube of **P** that protrudes into ℭ) about a vertical axis through its center would allow **T** to have been located in any of the other three cells on the top level. With ℭ containing the singletons **S** and **P**, it is doubly deficient and therefore can be tromino tiled by hypothesis. Lemma 2 shows that the remaining 2×2×3 boxes in the slab can be tiled, and as noted in the second paragraph of the proof of this Theorem, the rest of the shell can be tiled as well. Thus 𝔇 can be tiled.

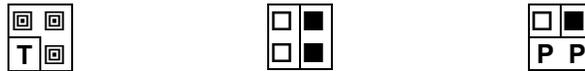

**Figure 17**

Figures 18 and 19 show the same types of views for the cases where **T** lies in the middle level of B and where **T** lies in the bottom level of B, respectively. Again, rotation of the box B and its protruding cube **P** lets **T** occupy any of the four positions on its level.

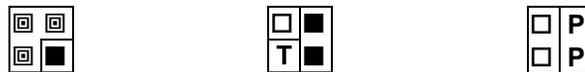

**Figure 18**

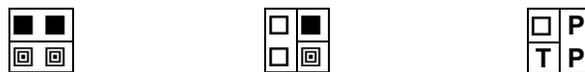

**Figure 19**

**Let T lie in a (3n+2)×(3n+2)×3 slab and 3n+2 be odd, with n > 1.**
Orient 𝔇 so the slab containing T is horizontal, on top of ℭ. Consider Figure 20, which shows in an overhead view a partition of this (3n+2)×(3n+2)×3 slab into a 3×3×3 box B in the upper right corner, 3×2×3 boxes in the row to the left of B, 2×3×3 boxes in the column below B, and 2×2×3 boxes in the remaining (3n−1)×(3n−1)×3 portion of the slab. By Lemmas 1 and 2, we know that the 3×3×3, 3×2×3, and 2×2×3 boxes can be tiled with trominoes. Suppose that T lies within one of the 2×2×3 boxes. We showed in the above analysis of a (3n+2)×(3n+2)×3 slab with 3n+2 even, that this 2×2×3 box can be tiled with three trominoes, T, and two cubes of a tromino whose third cube protrudes down into ℭ. In case T does not lie within the lower, left (3n−1)×(3n−1)×3 portion of the slab shown in Figure 20, a suitable rotation about a vertical axis will move this portion of the slab so it includes T. (This might not work if 3n+2 = 5, a case we treat separately below.) With the slab tiled and ℭ now doubly deficient, it follows as in the previous analysis for the case 3n+2 even, that all of 𝔇 can be tiled.

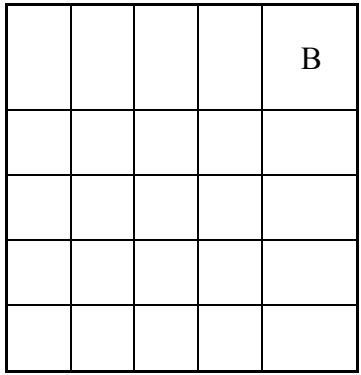

**Figure 20**

**Let T lie in a 5×5×3 slab.**
The left grid in Figure 21 shows an overhead view of the slab partitioned into a 3×3×3 cube, two 2×2×3 boxes, and four 2×1×3 boxes. If T is anywhere but above the center cell, it can (by rotation if necessary) be enclosed in a 2×2×3 box. Suppose T is in one of the three cells directly over the center of ℭ. Repartition the slab so that T lies in a 2×2×3 box as shown in the right grid in Figure 21. This figure also shows the tops of nine 1×2×3 boxes and one stack of three trominoes atop each other. In any case, the boxes not containing T can be tiled by Lemmas 1 and 2. As above, the 2×2×3 box containing T can be tiled with T, three trominoes, and two cubes of a tromino whose third cube protrudes down into ℭ, rendering ℭ doubly deficient and thus tilable. With the 5×5×3 slab and ℭ both tiled, all of 𝔇 can then be tiled.

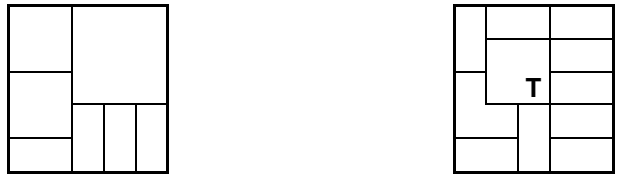

**Figure 21**

## Case B. T lies in a 3×3×(3n+2) rod or the 3×3×3 corner cube

For this one remaining case, note that the union of the 3×3×(3n+2) rod and the 3×3×3 cube amounts to a 3×3×(3n+5) rod. By Lemma 5, this extended rod is a union of 3×3×2 boxes (if 3n+5 is even) or one 3×3×3 cube together with 3×3×2 boxes. Moreover, because 3n+5 ≥ 8 the 3×3×3 cube may be repositioned, if necessary, so that **T** lies in a 3×3×2 box within the 3×3×(3n+5) rod.

Figure 22 shows an overhead view of the rod and with the adjacent 3×3×3 cube on its left end. Also shown is a view of the right face of a 3×3×2 box, with cells labeled for reference.

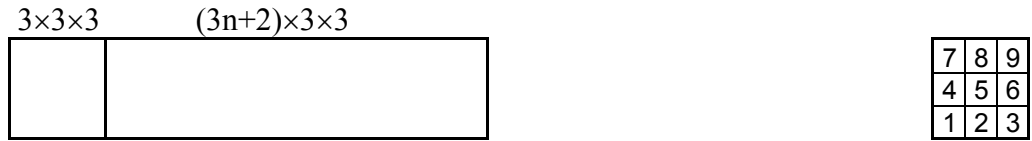

**Figure 22**

We now show that for **T** in the right 3×3 face of a 3×3×2 box B, it's possible to tile B with 5 trominoes, **T**, and a tromino having two cubes in B and its third cube protruding from the bottom level of B horizontally into the adjacent slab or, in case **T** is in the 3×3×3 corner cube, into another rod. If **T** lies, instead, in the left face, a similar tiling results from a simple reflection about a vertical plane, that is, about the hidden, middle 3×3 face of the 3×3×2 box.

We first tile B so as to allow **T** to lie in any of the four cells 2, 3, 5, and 6, by rotating the 2×2×1 array consisting of **T** and a tromino indicated by ▣'s in the first grid of Figure 23. This right side view shows **T** in cell 2, the tromino to be rotated, and parts of four other trominoes, one of which, marked **P**, protrudes outside the 3×3×2 box.

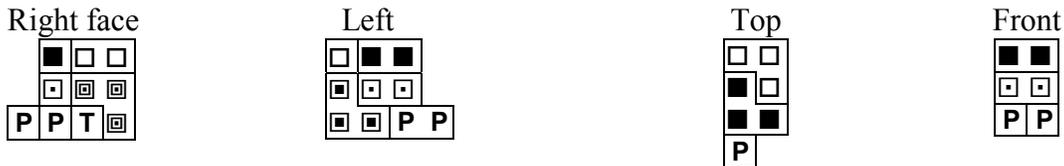

**Figure 23**

We next tile B so as to allow **T** to lie in either of cells 8 and 9, by rotating the 2×2×1 array on top, containing **T** and the ▣ tromino. The right side view in Figure 24 shows **T** in cell 9 and the overhead view shows how to carry out the rotation (clockwise, 90°.)

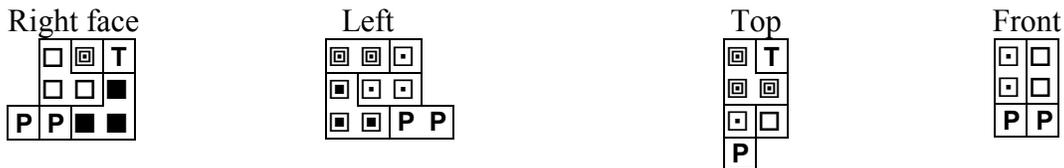

**Figure 24**

Figure 25 shows the tiling of B which allows **T** to lie in cell 7:

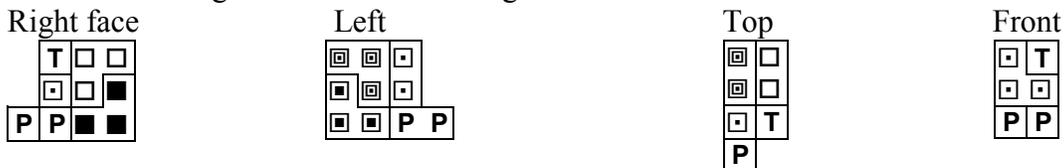

**Figure 25**

Finally, we tile B as in Figure 26 so that **T** may to lie in either of the two cells 1 and 4. The right side view shows the 2×2×1 cluster of **T** and 回, which rotates 90° clockwise to move **T** from cell 1 to cell 4.

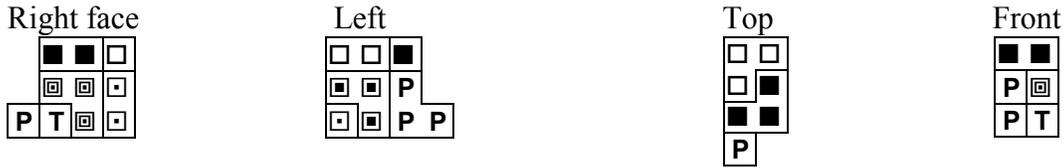

**Figure 26**

If **T** was in the 3×3×(3n+2) rod and *not* the 3×3×3 corner cube, we've noted that the union of the rod and corner cube consists of 3×3×2 boxes or else 3×3×2 boxes and one 3×3×3 cube, in such a way that **T** lies within one of the 3×3×2 boxes. We know how to tile the cube and those boxes not containing **T**, and we've just shown how to tile the 3×3×2 box containing **T** so as to leave a single **P** cube in an adjacent slab. In Case A we showed how to tile 𝔇 when a slab is singly deficient. Thus for **T** in the rod and not the corner cube, 𝔇 can be tiled.

The last case to be considered in this proof of Theorem 4 and hence of the General Theorem occurs when **T** lies in the 3×3×3 corner cube. Extend to the left the rod R1 of Figure 27, so it contains the 3×3×3 corner cube. Regard **T** as belonging to a 3×3×2 box in the extended R2. Tiling this box as just done above, the protruding **P** cube belongs to one of three positions in the end of the adjacent rod R2, as in the overhead view of Figure 27. We now place two trominoes which enable the tiling of 𝔇. One, indicated by 回's, will have a single cube on the bottom level of the adjacent (3n+2)×(3n+2)×3 slab (just atop ℭ), and two cubes occupying the bottom and middle level in rod R2 (shown as one 回 in the overhead views of Figure 27.) A second tromino will be placed so it has a single cube in the top level of ℭ (making ℭ doubly deficient and thus, by the induction hypothesis, tilable) and two cubes (indicated by ⊡) in the bottom level of the (3n+2)×(3n+2)×3 slab, just on top of ℭ. The three cubes (⊡,⊡, and 回) within the slab can easily be built up into a 3×2×2 box in the corner of the slab. As shown in Case A, when analyzing the (3n+2)×(3n+2)×3 slabs, such slabs can be partitioned so as to have a 3×2×2 box in any prescribed corner, with all the other boxes tiled as well.

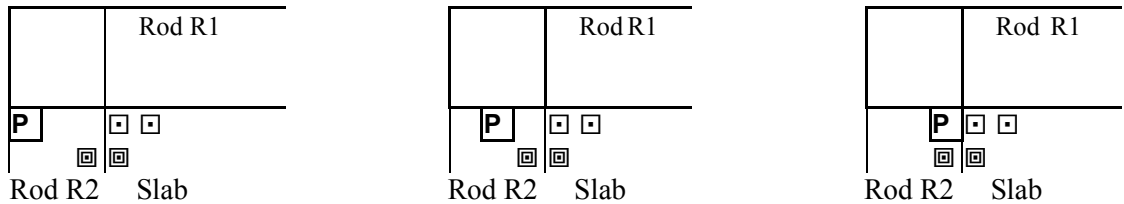

**Figure 27**

**P** and the two 回 cubes lie in a 3×3×2 box in rod R2 and having a face meeting a face of the 3×3×3 corner cube. As shown in our analysis of rods, R2 can be tiled so at one end (next to the 3×3×3 corner cube) there is a 3×3×2 box. For any of the three positions **P** shown in Figure 27, it's a simple exercise to complete the tiling (with five more trominoes) of the 3×3×2 box containing **P** and the two 回 cubes. Given that ℭ, the rod R2, the slab, and the rod R1 (containing the 3×3×3 corner cube) can be tiled when **T** lies in the corner cube, it remains only to tile the remaining rod and two remaining slabs. But Lemmas 3 and 4 assure this, so 𝔇 can be tiled, completing the proof of the General Theorem.

Acknowledgment: We thank Irene Starr ( www.starr.net/kbh/ ) for help with graphics and for logistic support.